\documentclass[a4paper,12pt]{article}
\usepackage{amssymb,delarray,amsmath}

\def\sf#1{{\mathsf{#1}}}
\let\3=\ss

\let\ss=\sss

\def\cc{{\mathbb C}}

\def\rr{{\mathbb R}}

\def\pp{{\mathbb P}}
\def\sph{{\mathbb S}}

\def\dim{\sf{dim}\,}

\def\lim{\mathop{\sf{lim}}}

\def\log{\sf{log}\,}

\def\reg{^\sf{reg}}
\def\sing{^\sf{sing}}

\def\sup{\sf{sup}\,}

\def\vol{\sf{Vol}}

\def\eps{\varepsilon}

\def\<{\langle}\let\la=\<
\def\>{\rangle}\let\ra=\>

\def\d{\partial}

\def\ddef{\mathrel{{=}\raise0.3pt\hbox{:}}}
\def\deff{\mathrel{\raise0.3pt\hbox{\rm:}{=}}}

\def\fraction#1/#2{\mathchoice{{\msmall{ #1\over#2}}}%
{{ #1\over #2 }}{{#1/#2}}{{#1/#2}}}
\def\norm#1{\left\Vert{#1}\right\Vert}

\def\le{\leqslant}

\def\emptyset{\varnothing}

\def\longpoints{\leaders\hbox to 0.5em{\hss.\hss}\hfill \hskip0pt}
\def\stateskip{\smallskip}
\def\state#1. {\stateskip\noindent{\bf#1. }} 
\def\statep#1. {\stateskip\noindent{\bf#1 }} 

\def\proof{\state Proof. \2}

\def\Chi{\raise 2pt\hbox{$\chi$}}

\def\ie{\hskip1pt plus1pt{\sl i.e.\/,\ \hskip1pt plus1pt}}

\def\sli{{\sl i)} } 
\def\slii{{\sl i$\!$i)} }

\def\sliii{{\sl i$\!$i$\!$i)} }

\def\reg{^\sf{reg}}
\def\sing{^\sf{sing}}

\def\v{\sf{v}}

\def\Chi{\raise 2pt\hbox{$\chi$}}

\let\phI=\phi\let\phi=\varphi\let\varphi=\phI
\def\calc{{\cal C}}

\def\calf{{\cal F}}
\def\calg{{\cal G}}

\def\calp{{\cal P}}

\def\calx{{\cal X}}

\def\eps{\varepsilon}

\def\d{\partial}

\def\1{{1\mkern-5mu{\rom l}}}

\def\ge{\geqslant}

\def\fraction#1/#2{\mathchoice{{\msmall{ #1\over#2}}}%
{{ #1\over #2 }}{{#1/#2}}{{#1/#2}}}

\def\le{\leqslant}

\def\emptyset{\varnothing}
\newcommand{\2}{\thinspace}

\newtheorem{thm}{Theorem}[section]
\newtheorem{lem}{Lemma}[section]
\newtheorem{prop}{Proposition}[section]
\newtheorem{corol}{Corollary}[section]
\newtheorem{defi}{Definition}[section]
\newtheorem{rema}{Remark}[section]
\newtheorem{quest}{Question}[section]

\def\eqqno(#1){\label{(#1)}}
\def\eqqref(#1){(\ref{(#1)})}
\def\qed{\ \ \hfill\hbox to .1pt{}\hfill\hbox to .1pt{}\hfill $\square$\par}

\advance\headheight 1.2pt

\def\el2{\sf{L^2}}

\title{Extra extension properties of equidimensional\\
holomorphic  mappings: results \\
and open questions}
\author{S. Ivashkovich}
\begin{document}
\maketitle

\begin{abstract}
Holomorphic (nondegenerate) mappings between complex manifolds of
the same dimension are of special interest. For example, they appear
as coverings of complex manifolds. At the same time they have very
strong ``extra`` extension properties in compare with mappings in
different dimensions. The aim of this paper is to put together the
known results on this subject, give some perspective on the general
strategy for future progress, prove some new results and formulate
open questions.
\end{abstract}

\section[1]{Introduction.}

\paragraph[1.1]{1.1. Results.} Most frequently a domain $D$ which admits a group 
of holomorphic authomorphisms 
$\Gamma$ acting on $D$ properly discontinuously withought fixed points is itself homogeneous,
\ie the group $Aut(D)$ of all biholomorphic authomorphisms of $D$ is transitive. 
In that case the quotient $X=D/\Gamma$ is a compact homogeneous manifold. It turns out that 
much more can be said about the extension of locally biholomorphic mappings 
in the case when $X$ is homogeneous and K\"ahler.

\smallskip Recall that a domain $(D,\pi)$ over a complex manifold $M$ os called 
locally pseudoconvex (locally Stein) if for every point $a\in \overline{\pi (D)}$
there exists a neighborhood $V\ni a$ such that all connected components of 
$\pi^{-1}(V)$ are Stein.

\bigskip\noindent\bf 
Theorem 1. {\it  Let $X$ be a locally homogeneous complex manifold and 
let $f:\hat D\to X$ be a meromorphic mapping from a locally pseudoconvex domain $(\hat D,\pi)$ 
over a complex manifold $M$. Suppose that $f$ is locally biholomorphic
outside of its indeterminacy set. Then $f$ is holomorphic (and therefore locally
biholomorphic) everywhere.

In, particular, if $f$ is a locally biholomorphic mapping from a domain $(D,\pi)$ over a Stein 
manifold $M$ to a compact locally homogeneous K\"ahler manifold then $f$ extends locally 
biholomorphically onto the envelpe of holomorphy $\hat D$ of $D$.}
\rm

\medskip\noindent The proof is given in Theorem 2.4 and in the Remark 2.4 after the proof of Theorem 2.4 
in  Section 2.
\begin{rema}\rm
The condition of compacity in this theorem, as in almost all results of this paper, can be relaxed to
disk-convexity. See more about this in Section 6.
\end{rema}
\begin{rema}\rm
For more results of this type see Corollary 3.2, Corollary 4.1, Corollaries 5.1, 5.2, 5.3 and Proposition 5.1.
\end{rema}

Let's now formulate an another type of result.

\bigskip\noindent\bf 
Theorem 2. {\it Let $p:X\to Y$ be a holomorphic fibration over a compact complex
manifold with compact K\"ahler fibers. Let $S\subset Y$ be a closed
subset such that $Y\setminus S$ is Stein. Then any meromorphic
section of our fibration, defined in a neighborhood of $S$
extends to a meromorphic section over the whole of $Y$.
}
\rm
\begin{rema}\rm
There is no assumption on how the K\"ahler metrics on fibers depend
on the point on the base. Of course the total space $X$ don't need
to be K\"ahler and even locally K\"ahler.
\end{rema}

\smallskip The proof, which is based on results of extension of meromorphic 
mappins into non-K\"ahler manifolds  is given in Section 4. More results are proved 
in Theorems 6.1 and  6.2.

\paragraph[1.2]{1.2. Open Questions.} 
A special accent in this paper is made to the open questions. This is done in order 
to present the authors vision of the future developpents in the area and in the hope
to get interested new people to enter the subject.

\smallskip Every section ends with few open questions relevant to thie section and 
the last Section 7 is entirely devoted to questions which concern the paper in 
general.

\begin{rema}\bf 1. \rm Mappings which decrease the dimension were studied in 
\cite{St} and \cite{Cha}.

\smallskip\noindent\bf 2. \rm Mappings which increase the dimension do not have any 
extension properties in general, see example in \cite{K1}. What one can expect at
the best is explained in Section 3.

\smallskip\noindent\bf 3. \rm We do not discuss here an extremely extended topic of 
locally biholomorophic (or proper) mappings between domains in $\cc^n$ and send the reader
to the recent paper \cite{NS} for the present state of art in this subject.
\end{rema}

\section[2]{Coverings of K\"ahler Manifolds}

\paragraph[2.1]{2.1. General K\"ahler Manifolds.}

We start with the relatively well understandable case of K\"ahler
manifolds. This case includes Stein manifolds, projective and
quasiprojective manifolds. Results quoted in this section are
directly applicable also to manifolds of class $\calc$, \ie
bimeromorphic to K\"ahler ones.

\smallskip In \cite{Iv1} the following theorem was proved:

\begin{thm}
Let $X$ be a  compact K\"ahler manifold. Then the following
conditions on $X$ are equivalent:

\smallskip
i) for any domain $D$ in a Stein manifold $M$, any holomorphic mapping $f : D
\longrightarrow X$ extend to a holomorphic mapping $\hat f : \hat D
\longrightarrow X$ from the envelope of holomorphy $\hat D$ of $D$
into $X$.

ii) $X$ doesn't contain rational curves, i.e., images 
of the Riemann sphere $\cc\pp^1$ under a non-constant holomorphic
mappings $\cc\pp^1\longrightarrow X$.
\end{thm}
\begin{rema}\rm
The condition on $X$ to be compact is too restricitive. It can be replaced 
by the disk-convexity, see section 6.
\end{rema}
As a corollary from this theorem we obtain a positive solution of the conjecture of
Carlson et Harvey, see \cite{CH}:

\begin{corol}
Let $D$ be a domain in a Stein manifold $M$ and let $\Gamma$ be a
subgroup of the group of holomorphic automorphisms of $D$ acting on
$D$ properly discontinuously without fixed points. If $D/\Gamma=X$ is
compact and K\"ahler then $D$ is Stein itself.
\end{corol}
Indeed, such $X$ cannot contain rational curves and therefore the
covering map extends from $D$ onto its envelope of holomorphy $\hat D$, 
which is Stein. Therefore
the only possibility is $\hat D=D$ and therefore $D$ is Stein itself. 
It can be viewed as certain
generalization of the theorem of Siegel, \cite{Sg}. In the theorem of
Siegel $D$ is supposed to be a bounded domain in  $M=\cc^n$ (as a
result in this case $X$ is, moreover, projective). If $D$ is not
necessarily bounded then $X$ may not be algebraic (example: a
non-algebraic torus as a quotient of $\cc^n$ by a lattice).

\smallskip\rm Recall the following 
\begin{defi}
A meromorphic mapping $f$ from a complex manifold $D$ to a complex manifold 
$X$ is an irreducible, locally irreducible analytic set $\Gamma_f\subset D\times X$ 
(graph of $f$) such that the natural projection $\pi_D:\Gamma_f\to D$ is proper and 
generically  one to one.
\end{defi} In that case there exists an analytic subset $I\subset D$ of codimension 
at least two such that $\Gamma_f \cap (D\setminus I)\times X$ is a graph of a 
holomorphic mapping (still denoted as $f$). This can be taken as a definition of 
a meromorphic mapping. 
The minimal $I$ satisfying this property is called the indeterminacy set of $f$.

In \cite{Iv2} the following conjecture of Griffiths,
see \cite{Gf}, was proved:

\begin{thm}
For any domain $D$ in a Stein manifold any meromorphic mapping from $D$
into a compact K\"ahler manifold $X$ extends to a meromorphic map
from the envelope of holomorphy $\hat D$ of $D$ into $X$.
\end{thm}
\begin{rema}\rm
Again, as in Theorem 2.1 one needs only disk-convexity from $X$ for
result to be still true.
\end{rema}
In the same way as above one obtains the following
\begin{corol}
Let $D$ be a domain in a complex manifold $M$ (not necessarily
Stein!) and $\Gamma$ a subgroup of the group of holomorphic automorphisms
of $D$ acting on $D$ properly discontinuously without fixed points.

(a) If $D/\Gamma$ is compact and K\"ahler, then $D$ is locally pseudoconvex.

(b) If $D/\Gamma=X'$ is Zariski  open in a compact K\"ahler $X$, then $D$
itself is Zariski open in a pseudoconvex domain $\hat D$ of $M$.
\end{corol}
If $D\subset \cc^n$
(b) is a theorem of Mok-Wong, \cite{MW}.

\medskip
\paragraph[2.2]{2.2. Homogeneous K\"ahler Manifolds.}

\begin{defi}
Complex manifold $X$ is called infinitesimally homogeneous if the
global sections of its tangent bundle generate the tangent space at
each point.
\end{defi}
One can prove, see Proposition 1.3 in [Hr], that for some natural
$N$ there exists a surjective endomorphism of holomorphic bundles
$\sigma : X\times \cc^N\to TX$. This property can be taken as a
definition of infinitesimally homogeneous manifold. All
parallelizable manifolds are inf. hom., as well as all Stein
manifolds and all complex homogeneous spaces under a real Lee group.
Every Riemann domain $(D,\pi)$ over an infinitesimally homogeneous
manifold is infinitesimally homogeneous itself.

\smallskip Using the morphism $\sigma$ and some riemannian metric on
$X$ one can define, as in [Hr] {\sl a  boundary distance } function
$d_D$ on $D$. Ruffly speaking $d_D(z)$ for $z\in D$ is the supremum
of the radii of balls $B$ with centers in $\pi (z)$ such that $\pi$
is injective over $B$. The principal result wee need from [Hr] is
contained in Theorem 2.1.
 It can be stated as follows:

\begin{thm}
If $(D,\pi)$ is a locally pseudoconvex domain with finite fibers over an
infinitesimally homogeneous complex manifold. Then the function
$-\log d_D(z))$ is plurisubharmonic.
\end{thm}
Note that no further assumptions on $X$ (like compactness or K\"ahlerness are needed).

\smallskip We shall need also the
Hironaka Resolution Singularities Theorem. We shall use the so
called embedded resolution of singularities, see \cite{Hi1}, \cite{BM}. Let us
recall  the notion of the sequence of blowings up over a complex
manifold $D$. Take a smooth closed submanifold $l_0\subset D_0:=D$
of codimension at least two. Denote by $\pi _1 : D_1\longrightarrow
D_0$ the blowing up of $D_0$ along $l_0$. Call this: \it a  blowing
up of $D_0$ along the closed center $l_0$. \rm  The exceptional
divisor $\pi ^{-1}(l_0)$ of this blowing up we denote by $E_1$.

We can repeat this procedure, taking a smooth closed submanifold
$l_1\subset E_1$ of codimension at least two in $D_1$ and produce $D_2$ and so on.

\begin{defi}
A finite sequence $\{ \pi ^j\} _{j=1}^N $ of such blowings up we
call {\sl a sequence of blowings up over $l_0\in D$},  or  {\sl a
regular modification over $l_0$}.
\end{defi}

\smallskip
By $\{ l_j\} _{j=0}^{N-1}$ we denote the corresponding centers and
by  $\{ E_j\}_{j=1}^N $ the exceptional divisors.
We put $\pi = \pi_1\circ ...\circ \pi_N$, $E$ denotes
the exceptional divisor of $\pi $, i.e. $E=\pi_N^{-1} (l_{N-1}\cup
...\cup (\pi_1\circ ... \circ \pi_N)^{-1}(l_0))$.

\smallskip Let $f:D\to X$ be a meromorphic mapping into a manifold $X$. Denote by $I$ the
set of points of indeterminacy of $f$,\ie $f$ is holomorphic on
$D\setminus I$ and for every point $a\in D$ $f$ is not holomorphic
in any neighborhood of $a$.

\smallskip\noindent\bf
Theorem. \it Let $f:D\to X$ be a meromorphic map between complex
manifolds $D$ and $X$. Then there exists a regular modification $\pi
: D_N\to D$ such that $f\circ\pi : D_N\to X$ is holomorphic.

\smallskip\noindent\rm See \cite{Hi1}. For the proof we refer also 
to \cite{BM}.

\begin{thm}
Let $X$ be a compact infinitesimally homogeneous K\"ahler manifold. Then every
locally biholomorphic mapping $f:D\to X$ from a domain $D$ over a
Stein manifold into $X$ extends to a locally biholomorphic mapping
$\hat f:\hat D\to X$ of the envelope of holomorphy $\hat D$ of $D$
into $X$.
\end{thm}
\proof Let $\hat f:\hat D\to X$ be the meromorphic extension of $f$.
Denote by  $I$ the set of points of indeterminacy of $\hat f$. Then
$\hat f|_{\hat D\setminus I}$ is locally biholomorphic and we can
consider the pair $(\hat D\setminus I, \hat f|_{\hat D\setminus I})$
as a Riemann domain over $X$.

This domain may not be locally pseudoconvex only at points of $I$.
But then its domain of existence $\tilde D$ over $X$
contains some part of the exceptional divisor $E$ of the
desingularization of $\hat f$.  The union of this part of $E$ with
$\hat D\setminus I$ is actually $\tilde D$ and the extension of
$\hat f_{\hat D\setminus I}$ to $\tilde D$ we denote as $\tilde f$. We consider $(\tilde
D, \tilde f)$ as a (locally  pseudoconvex) Riemann domain over $X$.

Suppose $\tilde D\setminus E$ is not empty. Then it is easy to
construct a sequence of analytic discs $\Delta_k$ in $\hat D\setminus I$ and
then in $\tilde D$ such that the boundaries of $\Delta_k$ stay in a
compact part of $\tilde D$, but $\Delta_k$ converge to a disc plus
some number of rational curves on $E\setminus \tilde D$. But this is
clearly forbidden by the plurisubharmonicity of $-\log d_{\tilde
D}$.

Therefore $(\tilde D, \tilde f)$ as a domain over $X$  coincides
with $(\hat D_N, \hat f_N)$ - desingularization of $\hat f$. But
then  $-\log d_{\tilde D}$ should be constant on all fibers of our
modification, because we can take as $\tilde D$ any locally pseudoconvex 
 neighborhoods of these fibers. This is impossible unless these fibers 
are poins. That means that $\hat f_N=\hat f$ and therefore $\hat f$ 
is holomorphic on $\hat D$.

\qed

\begin{rema}\rm
Theorem 2.2, Corollary 2.2 and Theorem 2.4 hold obviously  true for
manifolds of class $\calc$, \ie for manifolds that are bimeromorphic
to compact K\"ahler manifolds.
\end{rema}
\begin{rema}\rm
It is clear from the proof of Theorem 2.4 that the condition on $X$ to be compact can
be relaxed. In fact disk-convexity is sufficient, see Section 6. 
K\"ahlerness of $X$ was used also only once when we extended $f$ onto the
envelope of meromorphy. Therefore the Theorem 1 from the Introduction is 
also proved.
\end{rema}
\paragraph[2.3]{2.3. Open Questions.}

Let $X$ be a compact K\"ahler surface and let $f:B^*\to X$ be a locally biholomorphic
mapping of the punctured ball $B^*=B\setminus \{0\}$ into $X$. Then
$f$ extends meromorphically onto the whole ball $B$. The full image
by the extension $\hat f$ of the origin denote by $E:=\hat f[0]$.

\begin{quest}\rm
Prove that $E$ is an exceptional curve in $X$.
\end{quest}

\begin{quest}\rm What can be said about $\hat f[I]$ in the conditions 
of Theorem 2.4?
\end{quest}

\section[3]{Mappings into non-K\"ahler Manifolds}
\paragraph[3.1]{3.1. The Strategy.}
We start from the following remark. Let $X$ be a compact complex
manifold. Then due to the result of Gauduchon, see \cite{Ga}, $X$
admits a Hermitian metric $h$ such that its associated  form $\omega_h$
satisfies $dd^c\omega_h^{k}=0$, where $(k+1)$ is the complex dimension
of $X$.

\smallskip In fact we shall need a property which is easier to
prove:

\smallskip\noindent\it
Every compact complex manifold of dimension $k+1$ carries a strictly
positive $(k,k)$-form $\Omega^{k,k}$ with $dd^c\Omega^{k,k}=0$.

\smallskip\rm Indeed: either a compact complex manifold  carries a
$dd^c$-closed strictly positive  $(k,k)$-form or it  carries a
bidimension $(k+1,k+1)$-current $T$ with $dd^cT\ge 0$ but
$\not\equiv 0$. In the case of $\dim X=k+1$ such current is nothing
but a nonconstant plurisubharmonic  function, which doesn't exists
on compact $X$.

\smallskip
Let us introduce the class ${\cal G}_k$ of normal complex spaces,
carrying a nondegenerate positive $dd^c$-closed strictly positive
$(k,k)$-forms. Note that the sequence $\{ {\cal G}_k\}$ is rather
exhaustive: ${\cal G}_k$ contains all compact complex manifolds of
dimension $k+1$.

\smallskip Introduce furthermore the class of normal complex spaces
${\cal P}_k^{-}$ which carry a strictly positive $(k,k)$-form
$\Omega^{k,k}$ with $dd^c\Omega^{k,k} \le 0$. Note that ${\cal
P}_k^{-}\supset {\cal G}_k$. But Hopf three-fold $X^3=\cc^3\setminus
\{ 0\} /(z\sim 2z)$ belongs to ${\cal P}_1^{-}$ and not to ${\cal
G}_1$, see remark below.

\medskip Consider  the Hartogs figure
\begin{equation}
H^k_n(r):= \big[\Delta^n(1-r)\times \Delta^k \big] \cup\big[
\Delta^n \times A^k(r,1) \big]\subset \cc^{n+k}\,.
\end{equation}
Here $\Delta^{n}(r)$ stands for the $n$-dimensional polydisk of
radius $r$ and $A^k(r,1)=\Delta^{k}(1)\setminus \bar\Delta^{k}(r)$
for the $k$-dimensional annulus (or shell). In (1) one should think
about $0<r<1$ as being very close to $1$.

\medskip\noindent{\bf General Conjecture.} {\it  Meromorphic
mappings from $H_n^k(r)$ to compact (disk-convex) manifolds of class $\calp_k^-$ should
extend onto $\Delta^{n+k}\setminus A$ where $A$ is of Haussdorf
$(2n-1)$-dimensional measure zero. If the image manifold is from
class $\calg_k$ then $A\not= \emptyset$ should imply very restrictive
conditions on the topology and complex structure of $X$ (see results
below).
}
\paragraph[3.2]{3.2. Mappings into Manifolds of Class $\calg_1$.}
Let $A$ be a subset of $\Delta^{n+1}$ of Hausdorff
$(2n-1)$-dimensional measure zero. Take a point $a\in A$ and a
complex two-dimensional plane $P\ni a$ such that $P\cap A$ is of
zero length. A sphere $\ss^3= \{ x\in P : \Vert x-a\Vert =\eps \} $
with $\eps $ small will be called a "transversal sphere" if in
addition $\ss^3\cap A=\emptyset $.

\begin{thm}
Let $f:H_n^1(r)\to X$ be a meromorphic map into a compact complex
manifold $X$, which admits a Hermitian metric $h$, such that the
associated $(1,1)$-form $\omega_h$ is $dd^c$-closed (\ie $X\in\calg_1$). Then $f$
extends to a meromorphic map  $\hat f:\Delta^{n+1}\setminus A   \to
X$, where $A$ is a complete $(n-1)$-polar, closed subset of
$\Delta^{n+1}$ of Hausdorff $(2n-1)$-dimensional measure zero.
Moreover, if $A$ is the minimal closed subset such that $f$ extends
onto $\Delta^{n+1}\setminus A$ and $A\not= \emptyset $,  then for
every transversal sphere $\ss^3\subset \Delta^{n+1}\setminus A$  its
image $f(\ss^{3})$ is not homologous to zero in $X$.
\end{thm}

\begin{rema}
\bf 1. \rm A (two-dimensional) \it spherical shell \rm in a complex
manifold $X$ is the image $\Sigma $ of the standard sphere
$\ss^3\subset \cc^2$ under a holomorphic map of some neighborhood of
$\ss^3$ into $X$ such that $\Sigma $ is not homologous to zero in
$X$. Theorem 3.1 states that if the singularity set $A$ of our map
$f$ is non-empty, then $X$ contains spherical shells.
 
A good example to think about is a Hopf surface $H^2=\cc^2\setminus\{0\}/
(z\sim 2z)$ with the pluriclosed metric form $\omega = \frac{i}{2}
\frac{dz_1\wedge d\bar z_1 + dz_2\wedge\bar dz_2}{\norm{z}^2}$.

\smallskip\noindent\bf
2. \rm Consider now a Hopf three-fold $H^3=(\cc^3\setminus \{ 0\}
)/(z\sim 2z)$. The analogous metric form $\omega =\frac{i}{2}
\frac{dz_1\wedge d\bar z_1+ dz_2\wedge d\bar z_2+dz_3\wedge d\bar
z_3}{\Vert z\Vert^2}$ is not longer pluriclosed but only
plurinegative (i.e. $dd^c\omega \le 0$). Moreover, if we consider
$\omega $ as a bidimension $(2,2)$ current, then it will provide us
a natural obstruction for the existence of a  pluriclosed metric
form on $H^3$. That means that $H^3\in \calp^-_1\setminus \calg_1$. 
The natural projection $f:\cc^3\setminus \{ 0\} \to H^3$ has
singularity of codimension three and $H^3$ doesn't contains spherical
shells of dimension two (but contain a spherical shell of dimension
three). \cite{Iv4} contains extension theorem for mappings into
manifolds from class $\calp_1^-$ also.
\end{rema}

Later on in this paper we shall need one corollary from the
Theorem 3.1. A real two-form $\omega $ on a complex manifold $X$ is said
to "tame" the complex structure $J$ if for any non-zero tangent
vector $v\in TX$ we have $\omega (v,Jv)>0$. This is equivalent to
the property that the $(1,1)$-component $\omega^{1,1}$ of $\omega$
is strictly positive. Complex manifolds admitting a \it closed \rm
form, which tames the complex structure, are of special interest.
The class of such manifolds  contains all K\"ahler manifolds. On the
other hand, such metric forms are $dd^c$-closed. Indeed, if $\omega
=\omega^{2,0} +\omega^{1,1} +\bar \omega^{2,0}$ and $d\omega=0$,
then $\partial \omega^{1,1}=-\bar \partial \omega^{2,0}$. Therefore
$dd^c\omega^{1,1}=2i\partial \bar \partial \omega^{1,1}=0$. So the
Theorem 3.1 applies to meromorphic mappings into such manifolds. In
fact, the technique of the proof gives more:

\begin{corol}
Suppose that a compact complex manifold  $X$ admits
a strictly positive $(1,1)$-form, which is the $(1,1)$-component of
a closed form. Then every meromorphic map $f:H_n^1(r)\to X$ extends
onto $\Delta^{n+1}$.
\end{corol}

\begin{rema}1. \rm In particular, all results of Section 2 remain valid 
for such manifolds.

\smallskip\noindent\bf 2. \rm Theorem 3.1 stays valid for meromorphic mappings from 
all $H_n^k(r)$ for all $k\ge 1$. But it should be noted that in general 
extendibility of meromorphic mappings into some
complex manifold $X$ from $H_n^k(r)$ doesn't
imply extendibility of meromorphic mappings onto this $X$ neither from $H_{n+1}^{k}(r)$
no from $H_{n}^{k+1}(r)$
(for holomorphic mappings this is true), see example in  \cite{Iv6}.
\end{rema}

\paragraph[3.3]{3.3. Class $\calg_2$ and Dimension $3$.} The following result
was proved in \cite{IS}.

\begin{thm} Let $X$ be a compact complex space of dimension $3$
(more generally one can suppose that $X$ is of any dimension but carries a positive
$dd^c$-closed $(2,2)$-form). Then every meromorphic map
$f:H_1^2(r)\to X$ extends meromorphically onto $\Delta^{3}\setminus
A$, where $A$ is a zero-dimensional complete pluripolar set. If $A$
is non-empty then for every ball $B$ with center $a\in A$ such that
$\d B\cap A=\emptyset $, $f(\d B)$ is not homologous to zero in $X$,
\ie $f(\d B)$ is a spherical shell (of dimension $3$) in $X$.
\end{thm}
\begin{rema}\rm  Spherical shell of dimension $k$ in complex manifold
(space) $X$ is an image $\Sigma$ of the unit sphere
$\ss^{2k-1}\subset\cc^k$ under a meromorphic map $h$ from a
neighborhood of $\ss^{2k-1}$ into $X$ such that $\Sigma
=h(\ss^{2k-1})$ is not homologous to zero in $X$.
\end{rema}

\smallskip Results of such type have interesting applications to
coverings of compact complex manifolds as we shall see in the next
sections. From this theorem immediately follows that if
the covering manifold $\tilde V$ of a $3$-dimensional manifold $V$
is itself a subdomain in some compact complex manifold $Y$ then he
boundary of $\tilde V$ cannot have concave points.

\smallskip Let's give one more precise statement. Recall that a
complex manifold is called affine if it admits an atlas with affine
transition functions. In that case its universal covering is a
domain over $\cc^n$.

\begin{corol}
Let $V$ be a compact affine $3$-fold and let $(\tilde V, \pi)$ be
its universal covering considered as a domain over $\cc^3$ with
locally biholomorphic projection $\pi$. Then if $(\tilde V, \pi)$ is
pseudoconcave at some boundary point then $V$ contains a spherical
shell (of dimension $3$).
\end{corol}
Indeed, by the Theorem 3.2 the covering map $p:\tilde V\to V$ can be
extended to a neighborhood of a pseudoconcave boundary point, say
$a$, minus a zero dimensional set $A$. But this cannot happen unless
$\tilde V = V\cup A$ in a neighborhood of $a$. Therefore spheres
around $a$ project to shells in $V$ by the Theorem 3.2.

\begin{rema}\rm
Of course, an analogous result can be formulated for affine surfaces:
either the universal cover of an affine surface $V$ is Stein or
$V$ contains a spherical shell (of dimension two).
\end{rema}

\paragraph[3.4]{3.4. Open Questions.}
\begin{quest}\rm
We conjecture that the analogous result should hold for meromorphic
mappings in all dimensions. I.e. from $H_n^k(r)$ to compact
manifolds (and spaces)  in the classes ${\calp}_k^-$ and
${\calg}_k$. In particular, Theorem 3.2 should be true for
meromorphic mappings between equidimensional manifolds in all
dimensions.
\end{quest}
The main difficulty lies in the fact that it is impossible in
general to make the reductions (a)--(c) of \S 1 from \cite{IS}.
(Note that reductions (d)--(e) can be achieved in all dimensions.)

\begin{quest}\rm
One can start proving the general conjecture (as in Question 3.1) by considering
extension from $H_2^2(r)$ to a manifold of class $\calg_2$.
\end{quest}

\begin{quest}\rm
An analog of Corollary 3.2. in all dimensions seems to be an easier
problem then the Question 3.1 in its whole generality.
\end{quest}
It would be instructive to consult the paper \cite{BK} in this
regard.

\smallskip It is likely that one can say more about the singularity
set $A$ of the extended mapping in Theorems 3.1 and 3.2.
\begin{quest}\rm
Let $X$ is a compact complex manifold carrying a plurinegative
metric form, and let $f:\Delta^3\setminus S\to X$ is a meromorphic
mapping. Suppose that $A$ is a minimal closed subset of $\Delta^3$
such that $f$ extends onto $\Delta^3\setminus A$. Prove that each
connected component of $A$ is a complex curve.
\end{quest}
For general $X$ without special metrics the answer could be
negative, see examples in the last section of \cite{Iv4}.

\section[4]{Application to  K\"ahler Fibrations}

\paragraph[4.1]{4.1. Extension of Meromorphic Sections.} We start with
the proof of the Theorem 2 from the Introduction, which answers a question 
posed to the author by T. Ohsawa.

\proof \noindent{\sl Step 1.} {\it Every point $y\in Y$ has a
neighborhood $U$ such that $X_U=p^{-1}(U)$ - the union of fibers
over $U$, possesses a Hermitian metric such that its K\"ahler form
$\omega_U$ is a (1,1)-component of a closed form.}

To see this take a coordinate neighborhood $U\ni y$ such that $X_U$
is diffeomorphic to $U\times X_y$. Let $p_1:U\times X_y$ be the
projection onto the second factor. Let $\omega_y$ be a K\"ahler form
on $X_y$. Consider the following $1$-form on $X_U$: $\omega_U=
p^{*}dd^c|z|^2 + p_1^{*}\omega_y$, where $z$ is the vector of local
coordinates on $U$. $\omega_U$ is $d$-closed. Its $(1,1)$-component
is positive for $U$ small enough, since $\omega_y$ is positive on
$X_y$.

\smallskip Let $\rho$ be a strictly plurisubharmonic Morse exhaustion
function on the Stein manifold $W:=Y\setminus S$. Set $W_t = \{ y\in
W: \rho (y)> t\}$. Given a meromorphic section $\v$ on the
neighborhood of $S$. Then $\v$ is defined on some $W_t$. The set $T$
of $t$ such that $\v$ meromorphically extends onto $W_t$ is
non-empty and close.

\smallskip\noindent{\sl Step 2.} {\it $T$ is open.}

Let $t\in T$, then $\v$ is well defined and meromorphic on $W_t$.
Set $S_t=\{ y\in W:\rho (y)=t\}$. Fix a point $y_0\in S_t$. Take a
neighborhood $U$ of $y_0$ and form $\omega_U$ as in the {\sl Step
1}. If $y_0$ is a regular point of $S_t$ then there exists a Hartogs
figure $H\subset W_t$ such that the corresponding polydisk $D\ni
y_0$. By Corollary 3.1 the meromorphic mapping $\v:H\to D\times
X_{y_0}$ can be meromorphically extended to $D$ and we are done.

If $y_0$ is a critical point of $S_t$ then we use the result of
Eliashberg, see also Lemma 2.1 from \cite{FS}.  By this description
of critical points of strictly plurisubharmonic Morse functions we
can suppose that $y_0$ lies on a totally real disk $B$ in D and
$W_t\supset D\setminus B$. But then the argument remains the same,
because in this case we can also find an appropriate Hartogs figure
in $W_t$ with corresponding polydisk containing $y_0$.

Therefore $\v$ extends to $W_t$ for all $t$ and the Theorem is
proved.
\medskip

\qed

\begin{rema}\rm
Note that dealing with K\"ahler fibrations we were forced to use
Corollary 3.1 which concerns non-K\"ahler situation.
\end{rema}

\paragraph[4.2]{4.2. Non-K\"ahler deformations of K\"ahler manifolds.}
Recall that a complwx deformation of a compact complex manifold $X$
is a complex manifold $\calx$ together with a proper surjective
holomorphic map $\pi :\calx\to\Delta$ od rank one with connected
fibers and such that the fiber $X_0$ over zero is biholomorphic to
$X$. From \cite{Hi2} one knows that if $X_0$ is K\"ahler this
doesn't implies that the neighboring fibers $X_t$ are K\"ahler. But
the Step 1 in the proof of the Theorem 4.1 tells us that for $t\sim 0$
the fiber $X_t$ admits a Hermitian metric such that its associated
form is a $(1,1)$-component of a closed form. Therefore Corollary 3.1
applies to $X_t$.

\smallskip Let's give the formal statement. We say that a complex manifold 
$X$ possesses a meromorphic
extension property if for every domain $D$ in Stein manifold every
meromorphic mapping $f:D\to X$ meromorphically extends onto the
envelope of holomorphy of $D$.

\begin{corol}
Let $X_t$ be a complex deformation of a compact K\"ahler manifold $X_0$. 
Then for $t\sim 0$ $X_t$ possesses a meromorphic extension property.
\end{corol}

\paragraph[4.3]{4.3. Open Questions.}

\begin{quest}\rm
Suppose all $X_t$ for $t\not= 0$ possed meromorphic extension property (as, for
example, K\"ahler manifolds). Does $X_0$ possesses it as well?
\end{quest}

\begin{quest}\rm
If $X_0$ possesses a mer. ext. prop. does $X_t$ possesses it for $t$ close to zero?
\end{quest}

\section[5]{Coverings of non-K\"ahler Manifolds}
\paragraph[5.1]{5.1. General facts.} To stay within reasonable generality we shall restrict 
ourselves here with subdomains of $\cc\pp^n$ covering compact complex manifolds (this includes 
also subdomains of $\cc^n\subset\cc\pp^n$). Hovewer many statements have an obvious meaning
(reformulation) in the case of domains in general complex manifolds.

Locally pseudoconvex domains in (and over) both $\cc^n$ and $\cc\pp^n$ are Stein (with one
exception - $\cc\pp^n$ itself), see \cite{Ok,T}. They can cover both K\"ahler and non-K\"ahler 
manifolds. But Theorem 2.2 imply that:
\begin{corol}
If a  subdomain $D\subset\cc\pp^n$ covers a compact K\"ahler manifold $V$ then $D$ is Stein,
unless $V=\cc\pp^n$.
\end{corol}
An example of Stein domain covering a non-K\"ahler compact manifold is any Inoue surface
with $b_2=0$. Their universal covering is $\cc\times H$, where $H$ is the upper half-plane
of $\cc$.

\paragraph[5.2]{5.2. Coverings by domains from $\cc\pp^2$.} Since every compact complex surface admits a 
$dd^c$-closed metric  form, the Theorem 3.1 applies and we get:
\begin{corol}
If a subdomain $D\subset\cc\pp^2$ covers a compact complex surface $X$ then either $D$ is
Stein, or $\cc\pp^2$ itself, or $X$ contains a spherical shell.
\end{corol}
In $\cc\pp^3$ we have an analogous corollary from Theorem 3.2.
Recall that a domain $D\subset \cc\pp^n$ is $q$-convex if it admits
an exhaustion function such that its Levi form has at least $n-q+1$
strictly positive eigenvalues at each point outside of some compact
subset of $D$.
\begin{corol}
If $D\subset\cc\pp^3$ covers a compact complex $3$-fold then either $D$ is $2$-convex, or $D=\cc\pp^3$,
or $V$ contains a (three dimensional) spherical shell.
\end{corol}

\paragraph[5.3]{5.3. Coverings by ''large`` domains from $\cc\pp^3$.} A domain $D\subset \cc\pp^n$ is said to 
be ''large``
if its complement $\Lambda:=\cc\pp^n\setminus D$ is ''small`` in some sense. Different authors give different 
sense to the notion of being ''small``, see \cite{K2,L}  and therefore we shall reserve ourselves from
giving a general definition.

\smallskip Let's start from the remark that if $\Lambda\not=\emptyset$ then its Haussdorf $n$-dimensional 
(resp. $n-1$-dimensional) measure is non-zero if $n$ is even (resp. odd).  For example in $\cc\pp^2$ and in
$\cc\pp^3$ this condition is the same: $h_2(\Lambda)>0$, see \cite{L}. Both cases is easy to realise by 
examples. We have the following
\begin{prop}
Suppose a domain $D\subset\cc\pp^3$ covers a compact complex threefold $X$.

\smallskip\noindent{\sl Case 1.} If the complement $\Lambda =
\cc\pp^3\setminus D$ is locally a finite union of two-dimensional submanifolds, then $\Lambda$ is a 
union of finitely many lines.

\smallskip\noindent{\sl Case 2.} If the complement $\Lambda =
\cc\pp^3\setminus D$ is locally a finite union of three-dimensional submanifolds, then $\Lambda$ 
is foliated by lines.
\end{prop}
\proof Take a point $p$ on the limit set $\Lambda$ and find a point $q\in D$ and a sequence
of authomorphisms $\gamma_n\subset \Gamma$ such that $\gamma_n(q)\to p$. Here $\Gamma$
is a subgroup of $Aut(D)$ such that $D/\Gamma =X$. Due tp the Haussdorf dimension condition
on $\Lambda$ there exists a line $l\ni q$ such that $l\cap\Lambda =\emptyset$. Then
$\gamma_n(l)$ will converge to a line in $\Lambda$ passing through $p$.

\qed

In \cite{K2} an example of $\Lambda$ of dimension $3$ is constructed.

\paragraph[5.4]{5.4. Open Questions.}
\begin{quest}\rm
Prove an analog of the Case 1 of  Theorem 5.1 assuming only that $h_2(\Lambda)$ is finite.
\end{quest}

In that case the components of $\Lambda$ could be lines and points.

\begin{quest}\rm 
Suppose that the complement $\Lambda = \cc\pp^3\setminus D$ is locally a union of four-dimensional submanifolds. 
Are all components of $\Lambda$ necessariry either complex hypersurfaces or $CR$-manifolds of $CR$-dimension one?
Or, one can have components which are not $CR$-submanifolds? Are those $CR$-submanifolds Levi-flat?
\end{quest}

\cite{BK} contains an example on pp. 82-83 where one component of $\Lambda$ is a complex hyperplane,
 and another is  a Levi-flat ``perturbation'' of a complex hyperplane.

\section[6]{Disk-Convexity of Complex Spaces}
\paragraph[6.1]{6.1. The notion of disk-convexity.}
All results, except that of Section 4, presented in this paper are
valid for more general classes of complex manifolds and spaces then
just compact ones. Compacity can be replaced by much less restrictive
condition, namely by disk-convexity.
\begin{defi}
(a) Complex space $X$ is called disk-convex if for
every compact $K\subset X$ there is another compact $\hat K$ such
that  if for any holomorphic map $h:\bar\Delta \to X$ with
$h(\partial \Delta )\subset K$ one has $h(\Delta )\subset \hat K$.

(b) $X$ is called disk-convex in dimension $k$ if for
every compact $K\subset X$ there is another compact $\hat K$ such
that  if for any meromorphic map $h:\bar\Delta^k \to X$ with
$h(\partial \Delta^k)\subset K$ one has $h(\Delta^k)\subset \hat K$.
\end{defi}

\begin{rema}
\bf 1. \rm In all formulations of Section 2 ``compact K\"ahler'' can be replaced
by ``disk-convex K\"ahler''. Neither original proofs no backgrounds use
more then disk-convexity.

\smallskip\noindent\bf 2. \rm In formulations of Subsection 3.2 the same:
``compact of class $\calg_1$'' can be replaced by ``disk-convex of class $\calg_1$''.
This was actually done in \cite{Iv4}, see Theorem 2.2. there.

\smallskip\noindent\bf 3. \rm Theorem 3.2 is valid for  manifolds  from $\calg_2$
which are disk-convex in dimension $2$.
\end{rema}

\paragraph[6.2]{6.2. $k$-convexity $\Longrightarrow $  disk-convexity in dimension $k$.}

Now let us compare the notion of disk-convexity with other
convexities used in complex analysis. We shall see that our notion
is the most weaker one (and this is its great advantage).

\begin{defi}
A $\calc^2$-smooth real function $\rho $ on $X$ is called $k$-convex
if for any local chart $j:V \longrightarrow \tilde V\subset
\Delta^N$ there exists a real $\calc^2$-function $\tilde \rho $ on
$\Delta^N$ such that $\tilde \rho \circ j = \rho $ and the Levi form
of $\tilde \rho $ has at least $N-k+1$ positive eigenvalues at each
point of $\Delta^N$.
\end{defi}

\begin{defi} Complex space $X$ is called $k$-convex (in the
sense of Grauert) if there exists a $\calc^2$ exhaustion function
$\rho :X\longrightarrow [0,\infty [$, which is $k$-convex outside
some compact $K\Subset X$.
\end{defi}

\smallskip\rm
We shall start with the following

\smallskip\noindent\bf
Maximum Principle. {\it Let $\rho $ be a $k$-convex function on the
complex space $X$ and $A$ be a pure $k$-dimensional analytic subset
of $X$. If for some point $p\in A$ $\rho (p) = \sup_{a\in A}\rho
(a)$, then $\rho\mid_A\equiv const $.
}

\proof\rm If there is a smooth point $p\in A^{\reg}$ where
$\rho\mid_A$ achieves its maximum, then conclusion is clear. Really,
while the Levi form of $\rho_A:=\rho\mid_A$ has at least one
positive eigenvalue at $p$, one can find an analytically imbedded
disk $\Delta \ni p$ such that the restriction $\rho\mid_{\Delta }$
is subharmonic. This implies that  $\rho\mid_{\Delta }\equiv const$.
Further one can find a holomorphic coordinates
$(z_1,...,z_k)=(z_1;z^{'})$ in the neighborhood of $p$ such that
restriction $\rho_D$ of $\rho $ onto the every disk $D=\{
(z_1;z^{'}):z^{'}=0\} $ is subharmonic and such that our original
disk $\Delta $ is transversal to all such $D$. We conclude that
$\rho \equiv const$ in the neighborhood of $p$. The rest is obvious.

Now consider the case when $p\in A^{\sing}$ -- the set of singular
points of $A$. We shall be done if we shall prove that in the
neighborhood of $p$ there is another point $q\in A^{\reg}$  such
that $\rho (q)=\rho (p)$. Take a  neighborhood  $V$ of $p$ together
with imbedding $j:V\longrightarrow \tilde V\subset \Delta^N$ of $V$
as a closed analytic subset $\tilde V$ in the unit polydisk. Let
also $j(p)=0$. By $\tilde A$ let us denote $j(A\cap V)$ -- an
analytic subset of pure dimension $k$ in $\Delta^N$. Take some
irreducible component $B$ of $\tilde A$ passing through zero.

\begin{lem} Let $\Pi $ be a linear subspace of dimension $N-k+1$ of
$\cc^N$. Then for a subspace which is  a generic perturbation  of
$\Pi $ (and is again denoted as $\Pi $) there exists an $\eps >0$
such that $\tilde \Pi \cap B\cap \Delta^N_{\eps }$ is a complex
curve.
\end{lem}

\proof Blow up the origin in $\cc^N$. Let $\pp^{N-1}$ is an
exceptional divisor and $\pi :\cc^N\setminus \{ 0\} \to \pp^{N_1}$ a
natural projection. Denote by $\hat B$ and $\hat \Pi$ strict
transforms of $B$ and $\Pi $. Recall that $\pi^{-1}(\hat B\cap
\pp^{N-1})\cup \{ 0\} $ is a tangent  cone to $B$ at zero. While
$\hat B\cap \pp^{N-1}$ is of dimension $k-1$ and $\hat \Pi \cap
\pp^{N-1}$ is a linear subspace of dimension $N-k$, then for a
generic perturbation $\Pi $ the intersection $\hat\Pi\cap
\pp^{N-1}\cap \hat B$ is zerodimensional.

The usual properties of tangent cone imply  that $\Pi \cap B$ has
the tangent cone at zero of dimension one. And this implies that for
a small  enough $\eps >0$ the intersection this in $\Pi \cap B \cap
\Delta^n_{\eps } $ is a curve.

\smallskip
\hfill{Lemma is proved.}

\smallskip
Let us finish the proof of the maximum principle. While the Levi
form of $\tilde\rho :=\rho \circ j$ has at least $N-k+1$ positive
eigenvalues at zero, one can find a linear subspace $\Pi $ in
$\cc^n$ of dimension $N-k+1$ lying inside the positive cone of
${\cal L}_{\tilde \rho }(0)$. We can take instead of $\tilde A$ some
of its irreducible component $B$ passing through zero.  After a
small perturbation  $\Pi $ became transversal to $B^{\sing}$ still
being in the positive cone. Thus $\Pi \cap B^{\reg}\cap
\Delta^N_{\eps }\not= \emptyset $ for all $\eps >0$ small enough and
the same is true for the small perturbations of $\Pi $. Now our
lemma provides us with a perturbation $\Pi $ such that:

\noindent 1) $\Pi \cap B \cap \Delta^N_{\eps } =:C $  is a curve,
passing through zero for some $\eps >0$;

\noindent 2) $\Pi $ lies in the positive cone of ${\cal L}_{\tilde
\rho }(0)$;

\noindent 3) $C\cap B^{\reg}\not=\emptyset $.

\smallskip But this means that $\tilde\rho\mid_C$ is subharmonic.
Having zero as maximum it is constant. Thus we have found smooth
points where $\rho $ takes its maximum.

\smallskip
\hfill{q.e.d.}

\begin{thm}
$k$-convexity $\Longrightarrow $ disk-convexity in dimension $k$.
\end{thm}
\proof Let $\rho $ be an exhaustion function   on $X$, which is
$k$-convex outside compact $P$. Put $a=\sup_{x\in K\cup P}\rho (x)$,
and put $\hat K=\{ x\in X:\rho (x) \le a \} $. Let $h:\bar\Delta^k
\to X$ be some meromorphic map with $h(\partial \Delta^k)\subset K$.
Would $h(\bar \Delta^k)$ be not contained in $\hat K$ then $h(\bar
\Delta^k)\setminus \hat K$ would be a nonempty pure $k$-dimensional
analytic subset in $X\setminus \hat K$.

This clearly contradicts the maximum principle.

\medskip\qed

\begin{rema}\rm This Theorem answers the question which was posed to the Author by
D. Barlet. It is well known that $k$-convexity is nearly weakest notion among
convexities used in complex analysis.
\end{rema}

\paragraph[6.3]{6.3. Filling ``holes'' in Complex Surfaces.}
How fare can be a complex manifold or space be from being disk-convex?
This seems to be a difficult question. Here we shall indicate an interesting
particular case of being non-disk-convex. For the technical reasons we shall
restrict ourselves to complex dimension two. $B^*=B\setminus\{0\}$ will stand for the
punctured ball in $\cc^2$.

\smallskip Let $X$ be a normal complex surface, \ie a normal complex space of complex 
dimension tao, which will be supposed to be reduced and
countable at infinity. Following \cite{AS} we give the following

\begin{defi}
We say that $X$ has a hole if there exists a meromorphic mapping
$f:B^*\to X$ such that $\lim_{z\to 0}f(z)=\emptyset$.
\end{defi}
\begin{rema}\rm
If $X$ has a ``hole'' then it is certainly not disk-convex.
\end{rema}

But this particular cause of non-disk-convexity can be repaired.
\begin{thm}
Let $X$ be a complex surface. Then there exists a complex surface
$\hat X$ and a meromorphic injection $i:X\to \hat X$ such that

\sli $i(X)$ is open and dense in $\hat X$;

\slii $\hat X$ has no holes.
\end{thm}
\begin{rema}\rm
This result was announced in \cite{Iv5}, here we shall give the
sketch of the proof which crucially uses results of Grauert about
complex equivalent relations, see \cite{Gr1,Gr2}..
\end{rema}
\proof Let a ``hole'' $f:B^*\to X$ be
given.If there is a curve $C\subset B^*$ contracted by $f$ to a point
$p\in X$ we can blow-up $X$ at $p$ and get a new surface and a new map
which is not contracting $C$ and which is still a ``hole''. Since, after shrinking
$B$,there can be only finitely many contracted curves we can suppose without
loss of generality that
$$
f \text{ is  not contracting any curves in } B^*. \eqqno(*)
$$
On $B^*$ we define the following equivalence relation $x~y$ if
$f(x)=f(y)$. This means that if one of these points, say $y$ is an
indeterminacy point of $f$ then $f(x)\in f[y]$. If both $x$ and $y$
are points of indeterminacy then we require that $f[x]=f[y]$. This
equivalence relation $R\subset B^*\times B^*$ is an analytic set in
$B^*\times B^*$. This follows from the fact that $f$ is a ``hole''.
Really, one cannot have an accumulation point of $R$ of the kind
$(a,0)$ of $(0,a)$ with $a\not=0$. Moreover $R$ is semiproper for
the same reason. Therefore  $R$ extends to $B\times B$ it is a
meromorphic equivalence relation there in the sense of \cite{Gr2}.
By the results of \cite{Gr1,Gr2} the quotient $Q=B/R$ is a normal
complex surface.

\smallskip Now we can attach $Q$ to $X$ by $f|_R$ - quotient map and get a new normal surface
with a ``hole'' filled in.

\smallskip Using Zorn lemma one constructs a maximal extension $\hat X$ of $X$ such that
$X$ is open and dense in $\hat X$ ($\hat X$ is not unique!). The
``filling in'' procedure above implies that this $\hat X$ should be
Disk-convex.

\medskip\qed

\paragraph[6.4]{6.4. Open Questions.}
One could try to improve the result of Theorem 6.1:
\begin{quest}\rm
Can every complex surface be imbedded as a subdomain into a
disk-convex complex surface?
\end{quest}

In some cases another notions of ``disk-convexity'' are needed:

\smallskip\noindent
\sli A complex manifold $X$  is said to be disk-convex if for any
compact $K\subset X$ there exists a compact $\hat K\subset X$ such
that for every Riemann surface with boundary $(R,\d R)$ and every
holomorphic mapping $\phi :R\to X$ continuous up to the boundary the
condition $\phi (\d R)\subset K$ imply $\phi (R)\subset \hat K$.

\smallskip\noindent\slii $X$ is called disk-convex if for any
convergent on $\d \Delta$ sequence $\{\phi_n:\bar\Delta\to X\}$
of analytic disks this sequence converge also on $\bar\Delta$.

\smallskip\noindent\sliii The same definition can be given with sequences
of Riemann surfaces instead of the disk.
\begin{quest}\rm
What is the relation between all these notions and that defined in
Definition 6.1? Are they equivalent?
\end{quest}
Of course there are some obvious implications.

\section[7]{Open Questions}

\begin{quest} \rm Let the complex manifold $D$ is defined as two-sheeted
cover of $\Delta^2\setminus \rr^2$, i.e. $D$ is a "nonschlicht" domain
over $\cc^2$. Does there exist a compact complex manifold $X$ and a
holomorphic (meromorphic) mapping $f:D\to X$ which separates points?
\end{quest}

Note that the results of this paper imply that such $X$ if exists
cannot possed a plurinegative metric form. Thus examples could occur
starting from ${\dim}X\ge 3$.

\smallskip
In the following problems the space $X$ is equipped with some
Hermitian metric form $\omega$. On the subsets of $\cc^n$ the metric is always
$dd^c\Vert z \Vert^2$.

\begin{quest} \rm Consider a class ${\cal J}_R$ of meromorphic mappings
$f:\Delta^k\to X$, $X$ being compact, such that

(a) $\Vert Df\Vert\ge R>0$. Here $\Vert Df\Vert $ denotes the norm
of the differential of $f$;

(b) ${\vol}(f(\Delta^k)\le C_1$ for all $f\in {\cal J}_R$.

Prove that there is a constant $C_2=C_2(X,R,C_1)$ such
that ${\vol}(\Gamma_{f})\le C_2$ for all $f\in {\cal J}_R$.
\end{quest}

To estimate the volume of the graph of $f$ one should estimate the integral

\[
\vol (\Gamma_f) = \int_{\Delta^k}(dd^c||z||^2 + f^*\omega)^k = \sum_{j=0}^k
\int_{\Delta^k}(dd^c||z||^2)^j\wedge (f^*\omega)^{k-j},
\]
were only the first integral $\int_{\Delta^k}(f^*\omega)^k={\vol}(f(\Delta^k))$
is bounded by the condition of the question.

\smallskip The following question is of the same nature.

\begin{quest}\rm 
Let $f:\Delta^k_*\to X$ be a meromorphic mapping from a
punctured polydisk into a compact complex space $X$. Suppose that
${\vol}f(\Delta^k_*)<\infty $. Prove that $f$ meromorphically extends
to zero.
\end{quest}

\begin{quest} \rm  Let $f:\Delta^{k+1}_*\to X\in {\cal G}_k$ be a meromorphic
map from punctured $(k+1)$-disk into a compact complex space from
class ${\cal G}_k$. Prove that ${\vol}(f(A^k(r,1))=
O(\log^{\frac{k}{k-1}}(\frac{1}{r}))$ provided $k\ge 2$. In particular for
equidimensional maps $f:\Delta^n_*\to X^n$ one always should have
${\vol}(f(A^n(r,1))= O(\log^{\frac{n}{n-1}}(\frac{1}{r}))$.
\end{quest}

For $n=1$ there are no bounds on the growth of a meromorphic function in the 
punctured disk.

\begin{quest}
 \rm Fix some $0<r<1$ and some constant $R$. Fix also a compact
complex space $X$. Consider the following class ${\cal F}_R$ of meromorphic
mappings from $f:\Delta^n\to X$:

(1) ${\vol}_{2n}(\Gamma_f\cap (A^n(r,1)\times X))\le R$;

(2) for every $k$-disk $\Delta^k_z=\{ z\} \times \Delta^k$ ( where $z\in \Delta^{n-k}$)
${\vol}_{2k}(\Gamma_{f_z}\cap A^k_z(r,1)\times X)\le R$.

\medskip Prove that for any constant $l$ there is a constant $A$ such that for
any $f\in {\cal F}_R$ satisfying ${\vol}_{2k}(\Gamma_{f_z})\le l$ for
all restrictions $f_z$ of $f$ onto the $k$-disks $\Delta^k_z$ one
has ${\vol}_{2n}(\Gamma_f)\le A$.

\smallskip\noindent
(b) Vice versa: for any constant $a$ there is a constant $L$ such
that for any $f\in {\cal F}_R$ such that ${\vol}_{2n}(\Gamma_f)\le
a$ one has ${\vol}_{2k}(\Gamma_{f_s})\le L$ for all $\Delta^k_z$.
\end{quest}

The following question is a variation of questions 4.1 and 4.2. 
\begin{quest}\rm 
Let $\calx = \{X_t\}$ be a deformation of compact complex surfaces. Suppose that
$X_t$ for $t\not= 0$ contain a global spherical shell. Does $X_0$ contain a GSS?
\end{quest}

\begin{quest}\rm
Let $\calf$ be some family of holomorphic (meromorphic) mappings from the unit
polydisk $\Delta^{n+1}$ to a compact K\"ahler manifold $X$ (or $X\in \calg_1$ more generally).
Suppose that $\calf$ is equicontinuous on the Hartogs figure $H_n^1(r)$.
Will $\calf$ be equicontinuous on $\Delta^{n+1}$?
\end{quest}

See more about this question in \cite{Iv3}.

\ifx\undefined\bysame
\newcommand{\bysame}{\leavevmode\hbox to3em{\hrulefill}\,}
\fi

\def\entry#1#2#3#4\par{\bibitem[#1]{#1}
{\textsc{#2 }}{\sl{#3} }#4\par\vskip2pt}

\bigskip\noindent
U.F.R. de Math\'ematiques, \\
        Universit\'e de Lille-1\\
        59655 Villeneuve d'Ascq, France.\\
         E-mail: {\it ivachkov@math.univ-lille1.fr }

\bigskip\noindent
         IAPMM Acad. Sci. Ukraine,\\
         Lviv, Naukova 3b,\\
         79601 Ukraine.

\end{document}